\documentclass[11pt]{article}

\usepackage{amssymb,latexsym,amsmath}

\usepackage{graphicx}

\begin{document}

\newcommand{\bfi}{\bfseries\itshape}

\makeatletter

\@addtoreset{figure}{section}

\def\thefigure{\thesection.\@arabic\c@figure}

\def\fps@figure{h, t}

\@addtoreset{table}{bsection}

\def\thetable{\thesection.\@arabic\c@table}

\def\fps@table{h, t}

\@addtoreset{equation}{section}

\def\theequation{\thesubsection.\arabic{equation}}

\makeatother

\newtheorem{thm}{Theorem}[section]

\newtheorem{prop}[thm]{Proposition}

\newtheorem{lema}[thm]{Lemma}

\newtheorem{cor}[thm]{Corollary}

\newtheorem{defi}[thm]{Definition}

\newtheorem{hypo}{Brittle fracture hypothesis}[section]

\newtheorem{ehypo}[thm]{Equilibrium hypothesis (EH)}
\newtheorem{sehypo}[thm]{Strong equilibrium hypothesis (SEH)}

\newtheorem{rk}[thm]{Remark}

\newtheorem{exempl}{Example}[section]

\newenvironment{exemplu}{\begin{exempl}  \em}{\hfill $\surd$

\end{exempl}}

\newcommand{\comment}[1]{\par\noindent{\raggedright\texttt{#1}

\par\marginpar{\textsc{Comment}}}}

\newcommand{\todo}[1]{\vspace{5 mm}\par \noindent \marginpar{\textsc{ToDo}}\framebox{\begin{minipage}[c]{0.95 \textwidth}

\tt #1 \end{minipage}}\vspace{5 mm}\par}

\newcommand{\ea}{\mbox{{\bf a}}}
\newcommand{\eu}{\mbox{{\bf u}}}
\newcommand{\ueu}{\underline{\eu}}
\newcommand{\ueo}{\overline{u}}
\newcommand{\oeu}{\overline{\eu}}
\newcommand{\ew}{\mbox{{\bf w}}}
\newcommand{\ef}{\mbox{{\bf f}}}
\newcommand{\eF}{\mbox{{\bf F}}}
\newcommand{\eC}{\mbox{{\bf C}}}
\newcommand{\en}{\mbox{{\bf n}}}
\newcommand{\eT}{\mbox{{\bf T}}}
\newcommand{\eL}{\mbox{{\bf L}}}
\newcommand{\eV}{\mbox{{\bf V}}}
\newcommand{\eU}{\mbox{{\bf U}}}
\newcommand{\ev}{\mbox{{\bf v}}}
\newcommand{\eve}{\mbox{{\bf e}}}
\newcommand{\uev}{\underline{\ev}}
\newcommand{\eY}{\mbox{{\bf Y}}}
\newcommand{\eK}{\mbox{{\bf K}}}
\newcommand{\eP}{\mbox{{\bf P}}}
\newcommand{\eS}{\mbox{{\bf S}}}
\newcommand{\eJ}{\mbox{{\bf J}}}
\newcommand{\eB}{\mbox{{\bf B}}}
\newcommand{\leb}{{\cal L}^{n}}
\newcommand{\eI}{{\cal I}}
\newcommand{\eE}{{\cal E}}
\newcommand{\hen}{{\cal H}^{n-1}}
\newcommand{\eBV}{\mbox{{\bf BV}}}
\newcommand{\eA}{\mbox{{\bf A}}}
\newcommand{\eSBV}{\mbox{{\bf SBV}}}
\newcommand{\eBD}{\mbox{{\bf BD}}}
\newcommand{\eSBD}{\mbox{{\bf SBD}}}
\newcommand{\ecs}{\mbox{{\bf X}}}
\newcommand{\eg}{\mbox{{\bf g}}}
\newcommand{\paromega}{\partial \Omega}
\newcommand{\gau}{\Gamma_{u}}
\newcommand{\gaf}{\Gamma_{f}}
\newcommand{\sig}{{\bf \sigma}}
\newcommand{\gac}{\Gamma_{\mbox{{\bf c}}}}
\newcommand{\deu}{\dot{\eu}}
\newcommand{\dueu}{\underline{\deu}}
\newcommand{\dev}{\dot{\ev}}
\newcommand{\duev}{\underline{\dev}}
\newcommand{\weak}{\rightharpoonup}
\newcommand{\weakdown}{\rightharpoondown}
\renewcommand{\contentsname}{ }

\title{Equilibrium and absolute minimal states of Mumford-Shah 
functionals and brittle fracture propagation}
\author{Marius Buliga\footnote{"Simion Stoilow" Institute of Mathematics of the Romanian Academy,
 PO BOX 1-764,014700 Bucharest, Romania, e-mail: Marius.Buliga@imar.ro }}
\date{ }


\maketitle

\begin{abstract}
By a combination of geometrical and configurational analysis  
we study the properties of absolute minimal and equilibrium states  of general 
Mumford-Shah functionals, with applications to models of quasistatic brittle 
fracture propagation.  The main results 
 concern the mathematical relations between physical quantities 
 as energy release rate and energy concentration for 3D  
 cracks with complex shapes, seen as outer measures living on the 
 crack edge. 
 \end{abstract}


\section{Introduction}
\indent

A new direction of research in brittle fracture mechanics  begins with the 
article of  Mumford \& Shah  \cite{MS} regarding the problem 
of image segmentation. This problem,  which consists in finding the set of edges of a picture 
and constructing 
a smoothed version of that  picture, it turns to be intimately related to the problem of brittle 
crack evolution. 
In the before mentioned article Mumford and Shah propose the
following variational approach to the problem of image segmentation: 
let $g:\Omega \subset\mathbb{R}^{2} \rightarrow [0,1]$ 
be the original picture, given as a distribution of grey levels (1 is white and 0 is black), 
let $u: \Omega \rightarrow R$ 
be the smoothed picture and $K$ be the set of edges. $K$ represents the set where $u$ has jumps, 
i.e. $u \in C^{1}(\Omega \setminus K,R)$. The pair formed by the smoothed picture $u$ and the set 
of edges $K$ minimizes then the functional: 
$$I(u,K) \ = \ \int_{\Omega} \alpha \ \mid \nabla u \mid^{2} \mbox{ d}x \ + \ 
\int_{\Omega} \beta \  \mid u-g \mid^{2} \mbox{ d}x \ + \ \gamma \mathcal{ H}^{1}(K) \ \ .$$
The parameter $\alpha$ controls the smoothness of the new picture $u$, $\beta$ controls the 
$L^{2}$ distance between the smoothed picture and the original one and $\gamma$ controls the 
total length of the edges given by this variational method. The authors remark 
that for $\beta 
= 0$ the functional $I$ might be useful for an energetic treatment of fracture mechanics. 

An energetic approach to fracture mechanics is naturally suited to explain brittle crack 
appearance under imposed boundary displacements. 
The idea is presented  in the followings. 

The state of a brittle body is  described by a pair 
displacement-crack. $(\eu,K)$ is such a pair if $K$ is a crack --- seen as a surface 
---  
which appears in the body and  $\eu$ is a displacement of the broken body under the imposed 
boundary displacement, i.e. $\eu$ is continuous in the exterior of the surface $K$ and $\eu$
equals the imposed displacement $\eu_{0}$ on the exterior boundary of the body. 

Let us suppose 
that the total energy of the body is a Mumford-Shah functional of the form: 
$$E(\eu,K) \ = \ \int_{\Omega} w(\nabla \eu) \mbox{ d} x \ + \ F(\eu_{0},K) \ \ .$$
The first term of the functional $E$ represents the elastic energy of the body with the
displacement $\eu$. The second term represents the energy consumed to produce the crack 
$K$ in the body, with the boundary displacement $\eu_{0}$   as parameter. 
Then the crack that appears is supposed to be the second term of the pair 
$(\eu,K)$ which minimizes the total energy $E$. 

After the rapid establishment of mathematical foundations, starting with 
De Giorgi, Ambrosio \cite{DGA}, Ambrosio \cite{A1}, \cite{A2}, 
the development of such models continues with Francfort, Marigo \cite{FMa}, 
\cite{FMa2}, Mielke \cite{mielke}, Dal Maso, Francfort, Toader, \cite{dft05}, 
Buliga \cite{bu}, \cite{bu3}, \cite{bu2}.

In this paper we   introduce and study  equilibrium  and absolute minimal states of Mumford-Shah
functionals, in relation with a general model of quasistatic brittle crack 
propagation. 

On the space of the states of a brittle body, which are admissible with respect to 
an imposed Dirichlet condition, we introduce a partial order relation. Namely the state 
$(\eu,K)$ is "smaller than" $(\ev,L)$ if $L \subset K$ and $E(\eu,K) \leq E(\ev, L)$. 
Equilibrium states for the Mumford-Shah energy $E$ are then minimal elements of this
partial order relation. Absolute minimal states are just minimizers of the energy $E$. 

Both equilibrium states and absolute minimal ones are good candidates for solutions 
of models for quasistatic brittle crack propagation. Usually such models, based on 
Mumford-Shah energies, take into consideration only absolute minimal states. However, 
it seems to me that equilibrium states are better, because it is 
physically sound to define a state of equilibrium $(\eu,K)$ of a brittle body as one 
with the property that its total energy $E(\eu,K)$ cannot be lowered by increasing 
the crack further. 

For this reason we study here properties of equilibrium and absolutely minimal states 
of general Mumford-Shah energies. This study culminates with an inequality between 
the energy release rate and elastic energy concentration, both defined as outer measures 
living on the edge of the crack. 
This result  generalizes for tri-dimensional cracks 
with complex geometries  what is known about 
brittle cracks with simple geometry in two dimensions. In the two dimensional case, for 
cracks with simple geometry, classical use of complex analysis lead us to an 
equality between the energy release rate and elastic energy concentration at the tip 
of the crack. We prove that for absolute minimal states (corresponding to cracks 
with complex geometry) such an equality still holds, but for general equilibrium states 
we only have an inequality. Roughly stated, such a difference in properties of 
equilibrium and absolute minimal states comes from the mathematical fact that the 
class of first variations around an equilibrium state is only a semigroups. 
 
This research  might be relevant for 3D brittle fracture criteria applied for cracks 
with complex geometries. Indeed, it is very difficult even to formulate 3D fracture
criteria, because  in three dimensions  a crack of arbitrary shape 
does not have a finite number of "crack tips" (as in 2D classical theory),  
but an "edge" which is a collection of piecewise smooth curves in the 3D space.

\paragraph{Aknowledgements.} The  author  received  partial support 
from the Romanian Ministry of Education and Research, through the grant 
CEEX06-11-12/2006.

\section{Notations}

Partial derivatives of a function $f$ with respect to coordinate 
$\displaystyle x_{j}$ are denoted by $\displaystyle f_{,j}$. We use 
the convention of summation over the repeating indices. The open ball 
with center $\displaystyle x \in \mathbb{R}^{n}$ and radius $r>0$ is 
denoted by $B(x,r)$.

We assume that the body under study has an open, bounded, with locally Lipschitz
boundary,  reference configuration 
$\Omega \subset \mathbb{R}^{n}$, with $n=1, 2$ or $3$.  In the paper we shall 
use Hausdorff measures $\displaystyle \mathcal{H}^{k}$ in $\displaystyle 
\mathbb{R}^{n}$. For example, if $n=3$ then $\displaystyle \mathcal{H}^{n}$ is 
the volume measure, $\displaystyle \mathcal{H}^{n-1}$ is the area measure, 
$\displaystyle \mathcal{H}^{n-2}$ is the length measure. If $n=2$ then 
$\displaystyle \mathcal{H}^{n}$ is the area measure, $\displaystyle
\mathcal{H}^{n-1}$ is the length measure, $\displaystyle \mathcal{H}^{n-2}$ 
is the counting measure.

\begin{defi}
A smooth diffeomorphism with compact support in $\Omega$ is a function  
$\phi: \Omega \rightarrow \Omega$ with the following properties: 
\begin{enumerate}
\item[i)] $\phi$ is bijective;
\item[ii)] $\phi$ and $\phi^{-1}$ are  $C^{\infty}$ functions;
\item[iii)] $\phi$ equals the identity map of $\Omega$ near the boundary 
$\paromega$: 
$$supp \ ( id_{\Omega} - \phi ) \ \subset \subset \Omega \ \ .$$ 
\end{enumerate}
The set of all  diffeomorphisms with compact support in $\Omega$ 
is denoted by $\mathcal{D}$ or $\mathcal{D}(\Omega)$. 
\label{defd}
\end{defi}

The set $\mathcal{D}(\Omega)$  it is obviously non void because it
contains at least the identity map $id_{\Omega}$. Remark also that it is a 
group with respect to function composition. 

For any 
 $C^{\infty}$ vector field $\eta$ on $\Omega$ there is an unique associated 
 one parameter flow, which is a function $\phi: I \times \Omega 
 \rightarrow \Omega$, where $I \subset \mathbb{R}$ is  an open 
interval around $0 \in \mathbb{R}$, with the properties:
\begin{enumerate}
\item[f1)] $\forall t \in I$ the function $\phi(t,\cdot) \  = \ \phi_{t}(\cdot)$
   satisfies  i) and ii) from definition \ref{defd}, 
\item[f2)] $\forall t, t' \in I$, if $t-t' \in I$ then we have $\phi_{t'} \circ \phi_{t}^{-1} \ = \ 
\phi_{t-t'}$ ,
\item[f3)] $\forall \ t \in I$ we have $\eta \ = \ \dot{\phi}_{t} \circ \phi_{t}^{-1}$, where 
$\dot{\phi}_{t}$ means the derivative of $t \mapsto \phi_{t}$. 
\end{enumerate}

The vector field $\eta = 0$ generates the constant  flow 
$\phi_{t} \ = \ id_{\Omega}$. If 
$\eta$ has compact support in $\Omega$ then the associated flow 
$t \mapsto \phi_{t}$ is a curve in $\mathcal{D}$. 

A crack set $K$ is a piecewise Lipschitz surface with a 
boundary.  This means that exists  bi-Lipschitz functions 
$(f_{\alpha})_{\alpha \in 1 ... M}$, each of them defined over 
a relatively open subset $D_{\alpha}$ of  $\mathbb{R}^{n-1}_{+} \ = \ 
\left\{ y \in \mathbb{R}^{n-1} \mbox{ : } y_{n-1} \geq 0 \right\}$, with ranges in $\mathbb{R}^{n}$, such that: 
$$K \ = \ \cup_{\alpha = 1}^{M} f_{\alpha}(D_{\alpha}) \ \ , $$
$$ if \ \ \alpha \not = \beta \ \ then  \ \ 
f_{\alpha}(D_{\alpha} \setminus \partial \mathbb{R}^{n-1}_{+} ) \cap 
f_{\beta}(D_{\beta} \setminus \partial \mathbb{R}^{n-1}_{+} ) \ = \ \emptyset \ \ . $$ 
The edge of the crack $K$ is defined by
$$d K \ = \ \cup_{\alpha = 1}^{M} f_{\alpha}(D_{\alpha} \cap \partial \mathbb{R}^{n-1}_{+}) \ \ .$$ 
We shall denote further by $B_{r}(d K)$ the tubular neighborhood of radius $r$ of $d K$, given by
the formula: 
$$B_{r}(d K) \ = \ \cup_{x \in dK}  B(x,r) \ \ .$$
We denote by $[f] = f^{+} - f^{-}$ the jump of the function $f$ over the surface $K$ with respect 
to the field of normals $\en$.

\section{Mumford-Shah type energies}

\begin{defi}
 We describe the state of a brittle body by a pair  
$(\ev,S)$.  
The crack is seen as a piecewise Lipschitz surface $S$ in the topological 
closure $\displaystyle \overline{\Omega}$ of the reference 
configuration $\Omega$ of the body and 
$\ev$ represents the displacement of the body from the reference 
configuration. The displacement  $\ev$ has to be 
 compatible with the crack , 
i.e. $\ev$ has the regularity $C^{1}$ outside  the surface $S$.

The space of states of the brittle body with reference configuration 
$\Omega$ is denoted by $Stat(\Omega)$. 
\end{defi}

The main hypothesis in models of brittle crack propagation based on
Mumford-Shah type energies is the following. 

\smallskip

\noindent\textbf{Brittle fracture hypothesis.}
\textit{ The total energy of the 
 body subject to the  boundary displacement $\eu_{0}$ depends only on the state of the body 
$(\ev,S)$ and it has the expression: 
\begin{equation}
E(\ev,S) \ = \ \int_{\Omega} w(\nabla \ev) \mbox{ d} x \ + \ 
F(S ; \eu_{0}) \ \ . 
\label{genenrg}
\end{equation}
The first term of this functional is the  elastic energy associated to the 
displacement $\ev$; the second term represents the energy needed to produce 
the crack $S$, with the boundary displacement $\eu_{0}$ as parameter. }

We suppose that the elastic energy potential $w$ is a smooth, non negative 
function. 

 The most simple form of the function $F$ is the Griffith type energy: 
$$F(S ; \eu_{0}) \ = \ Const. \ \cdot \ Area \ (S) \ \ , $$
that is the energy consumed to create the crack $S$ is proportional, 
through a material constant, to the area of $S$. 
 
One may consider expressions of the surface  energy  $F$,  different from 
(\ref{genenrg}), for example:  
$$F(\ev, S) \ = \ \int_{S} \phi(\ev^{+},\ev^{-},\en) \mbox{ d}s \ \ ,$$
where $\en$ is a field of normals over $S$ , $\ev^{+}$, $\ev^{-}$ are the 
lateral limits of $\ev$ on $S$ with respect to directions $\en$, respective 
$-\en$ and $\phi$ has the property:
$$\phi(\ev^{+},\ev^{-},\en) \ = \  \phi(\ev^{-}, \ev^{+}, -\en) \ \ \ .$$ 
The function $\phi$, depending on the displacement of the "lips" of the crack, 
is a potential for surface forces acting on the crack. The expression 
(\ref{genenrg}) does not lead to such forces.

In general we shall suppose that the function $F$ has the properties:
\begin{enumerate}
\item[h1)] is sub-additive: for any two crack sets $A$, $B$ we have   
$$F(A \cup B ; \eu_{0}) \ \leq \  F(A ; \eu_{0}) \ + \ 
 F( B ; \eu_{0}) \ \ , $$
\item[h2)] for any $x \in \Omega$ and $r > 0$, let us denote by $\displaystyle 
\delta^{x}_{r}$ the dilatation  of 
center $x$ and coefficient 
$r$: 
$$\delta^{x}_{r} (y) = x + r(y-x) \quad .$$
Then, there is a constant $C \geq 1$ such that 
for any $A \subset \Omega$ with $F(A;\eu_{0}) < + \infty$ we have:  
$$F(\delta^{x}_{r}(A) \cap \Omega; \eu_{0} ) \ \leq \  C r^{n-1} F(A;\eu_{0})
\quad . $$
\end{enumerate}

The particular case $\displaystyle F(A;\eu_{0}) = G \hen(A)$ satisfies these 
two assumptions. In general these assumptions are satisfied for functions 
$\displaystyle F(\cdot ; \eu_{0})$ which are measures absolutely continuous with
respect to the area measure $\displaystyle \hen$.

A weaker property than h2), is the property h3) below. We don't explain here why
h3) is weaker than h2), but remark that h3) is satisfied by the same class of 
examples given for h2). 

For any $A \subset \Omega$, let us denote by $B(A,r)$ the tubular neighborhood of A: 
$$B(A,r) \ = \ \cup_{x \in A} B(x,r) \ \ . $$
We shall suppose that  $F$ satisfies: 
\begin{enumerate}
\item[h3)] for any $A \subset \Omega$ such that $F(A;\eu_{0}) < + \infty$, we have 
$$\limsup_{r \rightarrow 0} \frac{F(\partial B(A,r) \ \cap \ \Omega ; \eu_{0})}{r} \ < \ + \infty 
\ \ .$$
\end{enumerate}

\section{The space of admissible states of a brittle body}

\begin{defi}
The class of admissible  states of a brittle body with respect to the 
crack $\displaystyle F$ and with respect to the imposed 
displacement $\displaystyle \eu_{0}$ is defined as the  collection of all 
states $\displaystyle (\ev, S)$ such that 
\begin{enumerate} 
\item[(a)] $\displaystyle \eu = \eu_{0}$ on $\paromega \setminus S$, 
\item[(b)] $\displaystyle F \subset S_{u}$. 
\end{enumerate}
This class of admissible states is denoted by $\displaystyle Adm(F, \eu_{0})$. 
\label{defadmissible}
\end{defi}

An admissible displacement $\eu$ is a 
function which has to be equal to the imposed displacement on the boundary 
of $\Omega$ (condition (a)). Any such function $\eu$ is  reasonably smooth 
in the set $\displaystyle \Omega \setminus S_{u}$ and the function $\eu$ 
is allowed to have jumps along the set  $\displaystyle 
S$. 
Physically 
the set $\displaystyle $ represents the collection of all cracks in the 
body under the displacement $\eu$. The condition (b) tells us that the 
collection of all cracks associated to  an admissible displacement $\eu$  
contains $\displaystyle F$, at least. 

For some states $(\eu, S)$, the  crack set $S$ may have parts lying on the 
boundary of $\Omega$, that is 
$\displaystyle S \cap \paromega$ is a surface with positive area. 
In such cases we think about $\displaystyle S \cap \paromega$ 
as a region where the body has been detached from the machine which imposed 
upon the body the displacement $\displaystyle \eu_{0}$.

In a weak sense the  whole space of states of a brittle body may be identified 
with the space of special functions with bounded deformation $\eSBD(\Omega)$,
see \cite{BCDM}. 
Indeed, to every displacement field $\eu$ which is a special function with 
bounded deformation we associate the state of the brittle body described by 
$\displaystyle (\eu, \overline{\eS}_{u})$, where generally for any set $A$ we 
denote by $\overline{A}$ the topological closure of $A$. (Note that, technically,  
the crack set $\displaystyle \overline{\eS}_{u}$ may not be a collection of 
surfaces with Lipschitz regularity.)

On the space of states of a brittle body   we introduce a partial order relation. The definition is 
connected to definition \ref{defadmissible} and the brittle fracture 
hypothesis. 

\begin{defi}
Let $(\eu, S), (\ev, L) \in \, Stat(\Omega)$ be two states of a brittle body 
with reference configuration $\Omega$.  If
\begin{enumerate}
\item[(a)] $\displaystyle S \subset L$, 
\item[(b)] $\eu = \ev$ on $\paromega \setminus L$, 
\item[(c)] $\displaystyle E(\ev,L) \leq E(\eu,S)$, 
\end{enumerate}
then we write $\displaystyle (\ev, L) \leq (\eu, S)$. This is a partial order 
relation. 
\label{deford}
\end{defi}

There are many pairs $(\eu, S), (\ev, L) \in \, Stat(\Omega)$ such 
that $(\ev, L) \leq (\eu, S)$ and $(\eu, S) \leq (\ev, L)$, but 
$\eu \not = \ev$. Nevertheless 
such pairs have the same total energy $E$, the same crack set $S=L$, and 
$\eu = \ev$ on $\paromega \setminus L$. 

For a given boundary displacement $\displaystyle \eu_{0}$ and for given 
initial crack set $K$, on the set of admissible states 
$\displaystyle Adm(\eu_{0}, K)$ we have the same partial order relation.

\begin{defi}
An element  $\displaystyle (\eu, S) \in Adm(\eu_{0}, K)$ is minimal with 
respect to the partial order relation $\leq$ if for any  
$\displaystyle (\ev, L) \in Adm(\eu_{0}, K)$ the relation 
$(\ev, L) \leq (\eu, S)$ implies $(eu, S) \leq (\ev, L)$.

The set of equilibrium states with respect to given crack $K$ and imposed 
boundary displacement $\displaystyle \eu_{0}$ is denoted by $\displaystyle 
Eq(\eu_{0}, K)$ ant it consists of all minimal elements of $\displaystyle 
Adm(\eu_{0}, K)$ with respect to the partial order relation $\leq$. 

An element  $\displaystyle (\eu, S) \in   Adm(\eu_{0}, K)$ with 
the property that for any  
$\displaystyle (\ev, L) \in   Adm(\eu_{0}, K)$ we have 
$\displaystyle E(\eu, S) \leq E(\ev, L)$, is called an absolute minimal state. 
The set of absolute minimal states is denoted by $\displaystyle Absmin(\eu_{0},
K)$. 
\label{defequil}
\end{defi}

The physical interpretation of equilibrium states is the following. An
equilibrium state $\displaystyle (\eu, S) \in Eq(\eu_{0}, K)$ is one such that 
any other state $\displaystyle (\ev, L) \in Adm(\eu_{0}, K)$, which is 
comparable to $(\eu,S)$ with respect to the relation $\leq$, has the property 
$(\eu, S) \leq (\ev, L)$. In other words, equilibrium states are those with the
property: the total energy $E$ cannot be made smaller by prolongating the crack
set $S$ or by modifying the displacement $\eu$ compatible with the crack set 
$S$ and imposed boundary displacement $\displaystyle \eu_{0}$. 

Absolute minimal states are just equilibrium states with minimal energy. 

\begin{rk}
There might exist several minimal elements of of $\displaystyle 
Adm(\eu_{0}, K)$, such that any two of them are not comparable with respect to 
the partial order relation $\leq$. 
\end{rk}

For given expressions of the functions $w$ and $F$, we formulate the following 

\smallskip

\noindent\textbf{Equilibrium hypothesis (EH).}
\textit{ For any piecewise $C^{1}$ imposed 
boundary displacement $\eu_{0}$  and any crack $K$ the set of equilibrium states 
$\displaystyle Eq(\eu_{0}, K)$ is not empty. }

\smallskip

Without supplementary hypothesis on the total energy $E$, the EH does not 
imply that the set of absolute minimal states $\displaystyle 
Absmin(\eu_{0}, K)$ is non empty. Therefore the following hypothesis is stronger
than EH.

\smallskip

\noindent\textbf{Strong equilibrium hypothesis (SEH).}
\textit{For any piecewise $C^{1}$ imposed 
boundary displacement $\eu_{0}$  and any crack $K$ the set of equilibrium states 
$\displaystyle Absmin(\eu_{0}, K)$ is not empty.  }

\smallskip

\section{Models of quasistatic evolution of brittle cracks}
\label{themodels}

We shall describe here two models of quasistatic brittle crack propagation, 
according to Francfort, Marigo \cite{FMa}, \cite{FMa2}, Mielke \cite{mielke}, section 7.6, or Buliga 
\cite{bu2}, \cite{bu3}. At a first sight 
the models seem to be identical, but subtle differences exist. Further, 
instead of referring to a particular different model, we shall write about a 
general model of brittle crack propagation based on energy functionals, as if 
there is only one, general model, with different variants, according to the
choice among axioms listed further. 
Whenever necessary, the exposition will contain  variants of statements 
or assumptions which specializes the general model to one of the actual models
in use.

As an input of the model we have an initial crack set $\displaystyle 
K \subset \overline{\Omega}$ and a  curve of imposed displacements  
$\displaystyle t \in [0,T] \mapsto \eu_{0}(t)$ on the boundary of $\Omega$, the 
initial configuration of the body. 

We like to think about the configuration 
$\Omega$ as being an open, bounded subset of $\displaystyle \mathbb{R}^{n}$, 
$n = 1, 2, 3$, with sufficiently regular boundary (that is: piecewise Lipschitz 
boundary). 

The initial crack set $K$ has the status of an initial 
condition. Thus,  we suppose that $\displaystyle \partial \left(\mathbb{R}^{n} \setminus \Omega 
\right) = \paromega$.  For the same configuration $\Omega$ we may consider 
any crack set $\displaystyle  K \subset \overline{\Omega}$ as an initial crack. 
The crack set $K$ may be empty.

\begin{rk}
Models suitable for the evolution of brittle cracks under applied forces would 
be of great interest. Present formulations of the models of brittle crack 
propagation allows only the introduction of conservative force fields, as it 
is done in \cite{mielke} or \cite{FMa2}. The reason is that models based on 
energy minimization cannot deal with arbitrary force fields. In the case of a
conservative force field it is enough to introduce the potential of the force 
field inside the expression of the total energy of the fractured body. Thus, in 
this particular case we do not have to change substantially 
the  formulation of the model presented here, but only to slightly modify 
the expression of the energy functional. 
\end{rk}

In order to simplify  the model presented here,  we suppose that no 
conservative force fields are imposed on $\Omega$ or parts of $\paromega$. 
In the models described in \cite{mielke} or \cite{FMa2} such forces may be 
imposed.

\begin{defi}
A solution of the model is a curve of states of the brittle body 
$\displaystyle  t \in [0,T] \mapsto (\eu(t), S_{t})$ such that: 
\begin{enumerate}
\item[(A1)] (initial condition) $\displaystyle K \subset S_{0}$, 
\item[(A2)] (boundary condition) for any $t \in [0,T]$ we have 
$\displaystyle \eu(t) = \eu_{0}(t)$ on $\paromega \setminus S_{t}$, 
\item[(A3)] (quasistatic evolution) for any $t \in [0,T]$ we have 
$\displaystyle (\eu(t), S_{t}) \in \, Eq(\eu_{0}(t), S_{t})$, 
\item[(A4)] (irreversible fracture process) for any $t \leq t'$ we have 
$\displaystyle S_{t} \subset S_{t'}$, 
\item[(A5)] (selection principle) for any $t \leq t'$ and for any state 
$(\ev, S_{t}) \in \, Adm(\eu_{0}(t'), S_{t})$ we have $\displaystyle 
E(\ev, S_{t}) \geq E(\eu(t'), S_{t'})$.
\end{enumerate}
\label{defbul}
\end{defi}

From definition \ref{defequil} we see that (A2) is just a part of (A3). The
axiom (A2) is present in the previous definition only for expository reasons. 

The selection principle (A5) enforces the irreversible fracture process axiom 
(A4). Indeed, we may have severe non-uniqueness of 
solutions of the model. The axiom (A5) selects among all solutions satisfying 
(A1), ..., (A4), the ones which are energetically economical. The crack set 
$\displaystyle S_{t}$ does not grow too fast, according to (A5). For imposed 
displacement $\displaystyle \eu_{0}(t')$, the body with 
crack set $\displaystyle S_{t'}$ is softer than the same body with the crack set
$\displaystyle S_{t}$, for any $t \leq t'$.

As presented in definition \ref{defbul}, the model  has been proposed in 
Buliga \cite{bu2}.

In the models described in \cite{mielke}, \cite{FMa}, \cite{FMa2} 
we don't need the  selection principle (A5)  and the axiom (A3) takes the 
stronger form: 
{\it 
\begin{enumerate}
\item[(A3')] (quasistatic evolution) for any $t \in [0,T]$ we have 
$\displaystyle (\eu(t), S_{t}) \in  Absmin(\eu_{0}(t), S_{t})$.
\end{enumerate}
}

 \section{The existence problem}
 
The existence of equilibrium, or absolutely minimal states clearly depends 
on the ellipticity properties of the elastic energy potential $w$ (as shown 
for example in \cite{A2}, \cite{BCDM} or \cite{FMa}). This is related to the 
existence of minimizers of the elastic energy functional, as shown by relation 
(\ref{ech}) further on.  Some form of ellipticity 
of the function $w$ is sufficient, but it is not clear if such conditions are 
also necessary. Much effort, especially of a mathematical nature, has been spent
on this problem. 

In this paper we are not concerned with the existence problem,
however. Our purpose is to find general properties of solutions of brittle 
fracture propagation models based on Mumford-Shah functionals.  These 
properties  do not depend on particular forms of the elastic energy potential 
$w$, but on the hypothesis made in the general model. As any
 other model, the one studied in this paper is better fitted to some physical 
 situations than others. If some property of solutions of this model are 
 incompatible with a particular physical case, then we must deduce that the 
 model is not fitted for this particular case (meaning that at least one of 
 the hypothesis of the model is not suitable to this physical case). We are 
 thus able to provide a complementary information to  the one provided by 
 the existence problem. See further the Conclusions section for more on the
 subject.

\section{Absolute minimal states versus equilibrium states}

The differences between the models come from the difference between equilibrium
states and absolute minimal states.

Absolute minimal states are equilibrium states, but not any equilibrium state 
is an absolute minimal state.

Let us denote by $(\eu,S)$ an equilibrium state of the body, with respect to 
the imposed displacement $\displaystyle \eu_{0}$ and initial crack set $K$.

Consider first the class of all admissible 
pairs $(\ev, S')$ such that $S = S$. We have, as an application of definition
\ref{defequil}, then: 
\begin{equation}
\int_{\Omega} w(\nabla \eu) \mbox{ d}x \ \leq \ \int_{\Omega} w(\nabla \ev) \mbox{ d}x \ \ \ 
\forall \ \ev, \ \ev \ = \ \eu_{0} \mbox{ on } \paromega \setminus K \ , \ \ev \in C^{1}(\Omega 
\setminus K) \ \ . 
\label{ech}
\end{equation}
Thus any equilibrium state minimizes the elastic energy functional (in the 
class of admissible pairs with the same associated crack set). A sufficient  
condition for the existence of such minimizers is the polyconvexity of the 
elastic energy potential $w$.

The elastic energy potential function 
$\displaystyle w: M^{n\times n}(\mathbb{R}) \rightarrow\mathbb{R}$ associates to
any strain  $\displaystyle \eF \in M^{n\times n}(\mathbb{R})$ (here $n=2$ or
$3$)  the real value $w(\eF) \in
\mathbb{R}$. If this function    is smooth enough then we can define the
(Cauchy) stress tensor as
coming from the elastic energy potential:
$$\sig (\eu) \ = \ \frac{\partial w (\eF) }{\partial \eF}(\nabla \eu) \ \ .$$ 
The variational inequality  (\ref{ech}) implies that in the sense of 
distributions we have: 
$$ div \ \sig( \eu) \ = \ 0 \ \  $$ 
and that on the crack set $S$ we have $$\sig(\eu)^{+}\en \ = 
\ \sig(\eu)^{-} \en \ =  0 \ \ ,$$
where the signs $+$ and $-$ denotes the lateral limits of $\sig(\eu)$ with 
respect to the field of normals $\en$.

\subsection{Configurational relations for absolute minimal states}
\label{configsect}

We can also make smooth variations of the  
pair $(\eu,S)$. Here appears the first difference between absolute minimal and 
equilibrium states. We suppose further that $S \setminus K \not = \emptyset$, 
in fact we suppose that $S \setminus K$ is a surface with positive area.

If $\displaystyle (\ev,L) \in Adm(\eu_{0}, K)$ is an admissible state 
 and $\phi \in \mathcal{D}$ is a diffeomorphism of $\Omega$ with compact support,
 such that $K \subset \phi(K)$, 
then $(\ev \circ \phi^{-1}, \phi (S))$ is admissible too. 

If $(\eu, S)$ is an absolute minimal state then, as an application of definition
\ref{defequil}, 
 we have: 
\begin{equation}
E(\eu,S) \ \leq \ E( \eu \circ \phi^{-1}, \phi(S)) \ \ \ \forall \phi \in 
{\cal D} \, , \, K \subset \phi(K)  
\quad . 
\label{echdif}
\end{equation}
We may use (\ref{echdif}) in order to derive a first variation equality. 

We shall restrict further to the group $\mathcal{D}(K)$ of diffeomorphisms 
$\phi \in \mathcal{D}$ such that 
$supp \, \left( \phi - id\right) \cap K = \emptyset$. Vector fields $\eta$ which
generate one-parameter flows in $\mathcal{D}(K)$ are those with the property 
$supp \, \eta \cap K = \emptyset$. Further we shall work only with such vector
fields.

We shall admit further that for any smooth vector field $\eta$ there exist 
the derivatives at $t=0$ of the functions: 
$$t  \ \ \mapsto \ \ \int_{\Omega} w(\nabla (\eu \circ \phi_{t}^{-1}) ) 
\mbox{ d}x \ \ , \ \ t \ \ \mapsto \ \ F(\phi_{t}(K) ; \eu_{0}) \ \ \ ,$$ 
where $\phi_{t}$ is the one parameter flow generated by the vector field 
$\eta$. The relation  (\ref{echdif}) implies then: 
\begin{equation}  
\frac{d}{d t}_{|t=0} \  F(\phi_{t}(S) ; \eu_{0}) \ = \ - \ \frac{d}{d t}_{|t=0} 
\ \int_{\Omega} w(\nabla (\eu \circ \phi_{t}^{-1}) ) \mbox{ d}x \ \ . 
\label{fvar}
\end{equation}
Let us compute the right hand side  of (\ref{fvar}). We have 
$$  - \ \frac{d}{d t}_{|t=0} 
\ \int_{\Omega} w(\nabla (\eu \circ \phi_{t}^{-1}) ) \mbox{ d}x \ = \ 
\ \int_{\Omega}\left\{ - w(\nabla \eu) \ div \ \eta \ + \ \sig(\eu)_{ij} (\nabla \eu)_{ik} (\nabla 
\eta)_{kj} \right\} \mbox{ d}x \ \ . $$

For any vector field $\eta$, 
let us define, for any $x \in S$, $\displaystyle \lambda (x) =  \eta(x) \cdot 
\en(x)$, $\displaystyle \eta^{T}(x)  =  \eta(x) - \lambda)(x) \en(x)$,  
where $\en $ is a fixed field of normals over $S$.

With these notations, and recalling that the divergence of the stress field
equals $0$,  we have:
$$- \ \frac{d}{d t}_{|t=0} 
\ \int_{\Omega} w(\nabla (\eu \circ \phi_{t}^{-1}) ) \mbox{ d}x \ = \  \int_{S}  [w(\nabla \eu)]  
\lambda \mbox{ d } \hen \  + 
 $$ 
\begin{equation} 
+ \ \lim_{r \rightarrow 0} \int_{\partial B_{r}(dS)}  
 \left\{ [w(\nabla \eu)]  \lambda \ - \ [\sig(\eu)_{ij} (\nabla \eu)_{ik}]  \eta_{k}
\en_{j} \right\} \mbox{ d } \hen \ \ . 
\label{aspect}
\end{equation}

\begin{defi}
We introduce three kind of variations in terms of a vector field $\eta$ which 
generates an one parameter flow $\displaystyle \phi_{t} \in \mathcal{D}(K)$: 
\begin{enumerate}
\item[(a)] (crack neutral variations) for $\eta \ = \ 0$ on $S$; in this case 
we have $\displaystyle \phi_{t}(S) \ = \ S$ for any $t$,  
\item[(b)] (crack normal variations) for $\eta \ = \lambda \en$ on $S \setminus K$, 
with $\lambda: S \rightarrow \mathbb{R}$ a scalar, smooth function, such that 
$\lambda(x) = 0$ for any $x \in K \cup dS$, 
\item[(c)] (crack tangential variations) for $\eta \cdot \en = 0 $ on $S$.
\end{enumerate}
\label{defvaria}
\end{defi} 
For the case (a) of crack neutral variations  the relation (\ref{aspect}) 
 gives no new information, when compared with (\ref{ech}).

In the case (b) of crack normal variations, the relation (\ref{aspect}) 
implies  
$$\frac{d}{d t}_{|t=0} \  F(\phi_{t}(K) ; \eu_{0}) \ = \ \int_{S}  
[w(\nabla \eu)]  \lambda \mbox{ d } \hen \ \ .$$
In the particular case $F(S; \eu_{0}) \ = \ \hen (S)$ we obtain: 
$$\int_{S} \left\{ [w(\nabla \eu)] \ + H \right\} \lambda \mbox{ d} \hen \ = \ 0 \ \ , $$
where $H \ = \ - div_{s} \en \ = \ - \ div \ \en \ + \en_{i,j} \en_{i} 
\en_{j} $ is the mean curvature 
of the surface $S$. Therefore we have 
\begin{equation}
[w(\nabla \eu)(x)] \ + H(x) = 0
\label{meancurv}
\end{equation}
for any $x \in S \setminus K$.

In the case (c) of crack tangential variations, the relation (\ref{aspect})
implies 
$$\frac{d}{d t}_{|t=0} \  F(\phi_{t}(S) ; \eu_{0}) \ = \ $$
\begin{equation} 
= \ \lim_{r \rightarrow 0} \int_{\partial B_{r}(dS)}  
 \left\{ [w(\nabla \eu)]  \lambda \ - \ [\sig(\eu)_{ij} (\nabla \eu)_{ik}]  \eta_{k}
\en_{j} \right\} \mbox{ d} \hen \  .
\label{ji}
\end{equation}
This last relation admits an well known interpretation, briefly explained in the
next subsection. 

\subsection{Absolute minimal states for $n=2$}
\label{neq2}

Let us consider the case $n=2$ and the function 
$$F(S;\eu_{0}) \ = \ G \ \mathcal{H}^{1}(S) \ \ , $$
where $\mathcal{H}^{1}$ is the one-dimensional Hausdorff measure, i.e. the 
length measure. Let us suppose, for simplicity,  that the initial 
 crack set $K$ is empty  and the crack set $S$  of the absolute minimal state 
 $(\eu,S)$ has only one edge, i.e. $dS \ = \ \left\{ x_{0} \right\}$. 
Let us choose a vector field $\eta$ with compact support in $\Omega$ such that 
$\eta$ is tangent to $S$.  The equality (\ref{ji}) becomes then 
$$ G \ \eta(x_{0}) \cdot \tau (x_{0}) \ = \  
\lim_{r \rightarrow 0} \int_{\partial B_{r}(x_{0})}  
 \left\{ [w(\nabla \eu)]  \eta \cdot \en  \ - \ [\sig(\eu)_{ij} (\nabla \eu)_{ik}]  \eta_{k}
\en_{j} \right\} \mbox{ d} \hen \  ,  $$
where $\tau (x)$ is the unitary tangent in $x \in K$ at $K$.    
 If we suppose moreover 
that the crack $S$ is straight near $x_{0}$, and the material coordinates are 
chosen such that near $x_{0}$ we have $\eta(x) \ = \ \tau(x) \ = \ (1,0)$, 
then the equality (\ref{ji}) takes the form:
\begin{equation} 
 G \ = \ \lim_{r \rightarrow 0} \int_{\partial B_{r}(x_{0})}  
 \left\{ [w(\nabla \eu)]  \en_{1}  \ - \ [\sig(\eu)_{ij} (\nabla \eu)_{i1}]  
\en_{j} \right\} \mbox{ d} \hen \  .  
\label{grifull}
\end{equation}
We recognize in the right term of (\ref{grifull})  the integral $J$ of Rice; 
therefore at the edge of the crack the integral $J$ has to be equal to the 
constant $G$, interpreted as the constant of Griffith. 

The equality (\ref{grifull}) tells us that at the edge of a crack set belonging 
to an absolute minimal state  the Griffith criterion is fulfilled with equality.

\subsection{Configurational inequalities}

For equilibrium states which are not absolute minimal states we obtain just an 
inequality, instead of the equality from relation (\ref{ji}). Also, for such 
equilibrium states  there is no relation like (\ref{meancurv}) 
between the mean curvature of the crack set and the jump of elastic energy 
potential. We explain this further. 

The reason lies in the fact that if $\displaystyle (\eu, S) \in Eq(\eu_{0}, K)$ 
is an equilibrium state with $S \setminus K$ having positive area, and $\phi \in
\mathcal{D}(K)$ is a diffeomorphism preserving the initial crack set $K$, then 
we don't generally have the relation (\ref{echdif}). 

Indeed, in order to be able to compare $(\eu, S)$ with $\displaystyle 
(\eu \circ \phi^{-1}, \phi(S))$, we have to impose $S \subset \phi(S)$. Only for
these diffeomorphisms $\phi \in \mathcal{D}(K)$ the relation (\ref{echdif}) is
true. The class of these diffeomorphisms is not a group, like $\mathcal{D}(K)$, 
but only a semigroup. Technically, this is the reason for having only an
inequality replacing (\ref{ji}), and for the disappearance of relation 
(\ref{meancurv}). 

There is a necessary condition on the edge $dS$ of the crack set $S$, in order 
to have a trivial vector field $\eta$ which generates a one parameter flow  
$\displaystyle \phi_{t} \in \mathcal{D}(K)$ with $\displaystyle 
S \subset \phi_{t}(S)$ for any $t \in [0, T]$ (with $T>0$ sufficiently small). 
This condition is $\displaystyle dS \setminus K \not = \emptyset$.

Thus, for $\displaystyle (\eu, S) \in Eq(\eu_{0}, K)$ with $S\setminus K$ with
positive area, and $\displaystyle dS \setminus K \not = \emptyset$, we have 
\begin{equation}
E(\eu,S) \ \leq \ E( \eu \circ \phi_{t}^{-1}, \phi_{t}(S)) \ \ \ \forall t \in 
[0,T]   \quad ,  
\label{inech}
\end{equation}
for any one parameter flow $\displaystyle \phi_{t} \in \mathcal{D}(K)$ with 
$\displaystyle S \subset \phi_{t}(S)$ for any $t \in [0, T]$. 

In relation (\ref{inech}) crack normal variations (case (b) of definition
\ref{defvaria}) are prohibited. But these type of variations led us to the
relation (\ref{meancurv}). We deduce that for an equilibrium state 
$\displaystyle (\eu, S) \in Eq(\eu_{0}, K)$ , such that $S\setminus K$ has 
positive area, and $\displaystyle dS \setminus K \not = \emptyset$, the 
relation (\ref{meancurv}) does not necessarily hold. 

The crack tangential variations (case (c) of definition \ref{defvaria}) are
allowed in relation (\ref{inech}) only for $t \geq 0$. That is why we get only 
a first variation inequality: 
$$\frac{d}{d t}_{|t=0} \  F(\phi_{t}(S) ; \eu_{0}) \ \geq \ $$
\begin{equation} 
\geq \ \lim_{r \rightarrow 0} \int_{\partial B_{r}(dK)}  
 \left\{ [w(\nabla \eu)]  \lambda \ - \ [\sig(\eu)_{ij} (\nabla \eu)_{ik}]  \eta_{k}
\en_{j} \right\} \mbox{ d} \hen \  , 
\label{jineq}
\end{equation}
for any vector field $\eta$ which generates one parameter flow 
$\displaystyle \phi_{t} \in \mathcal{D}(K)$ with 
$\displaystyle S \subset \phi_{t}(S)$ for any $t \in [0, T]$.

The physical interpretation of relation (\ref{jineq}) is the following: 
the crack set $S$ of an equilibrium state satisfies the Griffith criterion of 
fracture, but, in distinction with the case of an absolute minimal state, 
there is an inequality instead of the previous equality. We are aware of at
least one example where this inequality is strict. This case concerns a crack 
set in 3D formed by a pair of intersecting, transversal planar cracks. Such a
crack set has an edge (in form of a cross), but also a "tip" (at the intersection 
of the edges of the planar cracks.  The physical implications of the 
inequality (\ref{jineq}) are that such a 3D crack may propagate in different 
ways, either along a crack tangential variation, or along a more topologically 
complex shape, by loosing its "tip". An article in preparation is dedicated 
to this subject.

We may  interpret the Griffith criterion of fracture, in the 
form given by relation (\ref{jineq}),  as a first order stability 
condition for the crack $S$ associated to the state of a brittle body. 
Surprisingly then,  absolute minimal states are first order neutral (stable 
and unstable), even 
if globally stable (as global minima of the total energy). There might exist 
equilibrium states for which we have strict inequality in relation 
(\ref{jineq}). Such states are surely not absolute minimal, but they seem 
to be first order stable, if our interpretation of (\ref{jineq}) is 
physically sound.

\subsection{Concentration of energy from comparison with admissible states}

We can obtain energy concentration estimates   from comparison of the energy 
of the equilibrium state $\displaystyle (\eu,S) \in Eq(\eu_{0}, K)$ with other 
particular admissible pairs.

Let $x_{0} \in \Omega$ be a fixed point and $r > 0$ such that $\displaystyle 
B(x_{0}, r) \subset \Omega$. 
We construct the following admissible pair $(\ev_{r},S_{r})$:
$$\ev_{r} (x) \ = \ \left\{ \begin{array}{ll}
\eu (x) & \mbox{ if } x \in \Omega \setminus B(x_{0},r) \\
0 & \mbox{ if } x \in \Omega \cap B(x_{0},r) \ \ , 
\end{array}
\right. $$
$$S_{r} \ = \ S \cup \partial B(x_{0},r) \ \ . $$
We have then the inequality $E(\eu,S) \ \leq \ E(\ev_{r}, S_{r})$, for any 
$r > 0$ sufficiently small. 
We use the properties h1), h2) of  $F$ to deduce that for 
any $x_{0} \in \Omega$ and $r > 0$ we have : 
\begin{equation}
\int_{B(x_{0},r)} w(\nabla \eu) \mbox{ d} x \ \leq \ C
\Omega_{n}(x;\eu_{0})  \ r^{n-1} \quad ,   
\label{grp1}
\end{equation}  
where $\displaystyle \Omega_{n}(x_{0};\eu_{0})$ is a number defined by 
$$\Omega_{n}(x_{0};\eu_{0}) = F(\partial B(x_{0}, 1); \eu_{0}) \quad .$$
In the case  of Griffith type surface energy  $\displaystyle F(S; \eu_{0}) = 
G \hen(S)$ we have 
$$\Omega_{n}(x_{0};\eu_{0}) = G \omega_{n} \quad , $$
with $\displaystyle \omega_{n}$ the area of the boundary of the unit ball 
in $n$ dimensions, that is $\displaystyle \omega_{1} = 2$, $\displaystyle
\omega_{2} = 2 \pi$, $\displaystyle \omega_{3} = 4 \pi^{2}$.

This inequality lead us to the following energy concentration property for 
$\eu$: 
\begin{equation}
\limsup_{r \rightarrow 0} \frac{ \int_{B(x_{0},r)} w(\nabla \eu) \mbox{ d} x}{ r^{n-1}} 
\ \leq \ C \Omega_{n}(x_{0}; \eu_{0}) \quad .
\label{grp2}
\end{equation}

The term from the left hand side of 
the relation (\ref{grp2})  is the concentration factor of the elastic energy 
around the point $\displaystyle x_{0}$.

The relation (\ref{grp2}) shows that the distribution of elastic energy of the 
body in the state $(\eu, S)$ is what we expect it to be, from the physical
viewpoint. Indeed, let us go back to the case $n=2$.  
It is well known that in the case of linear elasticity in two dimensions, 
if $(\ev,S)$ is a pair displacement-crack such that $div \ \sig (\ev) \ = \ 0$ 
outside $S$ and $\sig (\ev)^{+} \en  \ = \ \sig (\ev)^{-} \en \ = \ 0 $ 
on $S$ then $\ev$ behaves like $\sqrt{r}$ near the edge of the crack, hence 
the elastic energy behaves like  $\displaystyle r^{-1}$. We recover then the 
relation (\ref{grp2}) for $n=2$.

The relation (\ref{grp2}) does  imply that elastic energy concentration has an upper bound, 
but it does not imply that the energy concentration is positive at the tip of the crack. 
In the case $n=2$, for example, and for general form of the elastic energy density, 
 the relation (\ref{grp2}) tells us that if there is a concentration of energy 
 (that is if the density of elastic energy goes to infinity around the point $x$ in the 
 reference configuration) then the elastic energy density behaves like $\displaystyle r^{-1}$. 
 But it might happen that the elastic energy density is nowhere infinite. In this case 
 we simply have 
 $$\limsup_{r \rightarrow 0} \frac{ \int_{B(x_{0},r)} w(\nabla \eu) \mbox{ d} x}{ r^{n-1}} \ = \ 0$$ 
 which is not in contradiction with (\ref{grp2}).

From the hypothesis h3) upon the surface energy $F$ we get a slightly different 
estimate. We need first a definition.

\begin{defi}
For the equilibrium state $\displaystyle (\eu,S) \in Eq(\eu_{0}, K)$  and for 
any open set 
$A \subset \Omega$ we define:
$$CE(\eu,S)(A) \ = \ \limsup_{r \rightarrow 0}  
\frac{ \int_{B((dS \cap A ,r) \cap \Omega} w(\nabla \eu) \mbox{ d} x}
{r} \ \ , $$
$$CF(S; \eu_{0}) ( A) \ = \ \limsup_{r \rightarrow 0} \frac{ F(\partial B(dS 
\cap A,r); \eu_{0})}{r}
 \ \ . $$
 
The functions $\displaystyle CE(\eu,S)(\cdot)$, $\displaystyle
CF(S;\eu_{0})(\cdot)$ are  sub-additive functions which by well-known 
techniques induce  outer measures over the $\sigma$-algebra of 
borelian sets in $\Omega$. 

The function $\displaystyle CE(\eu,S)(\cdot)$ is called the elastic energy
concentration measure associated to the equilibrium state $(\eu,S)$. Likewise, 
the function  $\displaystyle
CF(S;\eu_{0})(\cdot)$ is called the surface energy concentration measure
associated to $(\eu,S)$. 
\label{deficonc}
\end{defi}

\begin{thm}
Let $\displaystyle (\eu, S) \in Eq(\eu_{0},K)$ be an equilibrium state. Then 
for any open set $A \subset \Omega$ we have 
$$CE(\eu,S)(A) \leq CF(S;\eu_{0})(A) \quad .$$
\label{firsthm}
\end{thm}

\paragraph{Proof.}
We consider, for any closed subset $A$ of $\Omega$ the following admissible 
state $\displaystyle (\eu_{r, A}, S_{r,A})$ given by: 
 $$\eu_{r, A} (x) \ = \ \left\{ \begin{array}{ll}
\eu (x) & \mbox{ if } x \in \Omega \setminus B(dS \cap A ,r) \\
0 & \mbox{ if } x \in \Omega \cap B(dS \cap A,r) \ \ , 
\end{array}
\right. $$
$$S_{r, A} \ = \ S  \cup \partial B(dS \cap A,r) \ \ . $$
The state $(\eu,S)$ is an equilibrium state and $\displaystyle (\eu_{r, A},
S_{r,A})$ is a comparable state, therefore we obtain: 
$$\int_{B(dS \cap A ,r) \cap \Omega} w(\nabla \eu) \mbox{ d} x \ \leq \ 
F(\partial B(dS \cap A,r); \eu_{0}) \quad .$$ 
We get eventually: 
$$ \limsup_{r \rightarrow 0}  \frac{ \int_{B(dS \cap A ,r) \cap \Omega} 
w(\nabla \eu) \mbox{ d} x}
{r}  \ \leq \ \limsup_{r \rightarrow 0} \frac{ F(\partial B(dS \cap A,r); 
\eu_{0})}{r} \quad . \quad \quad \square $$

Theorem \ref{firsthm}  shows that an equilibrium state satisfies a kind of Irwin type 
criterion. Indeed, Irwin criterion is formulated in terms of  stress
intensity factors. Closer inspection reveals that really  it is 
formulated in terms of elastic energy concentration factor, and that for 
special geometries of the crack set, and for linear elastic materials, 
we are able to compute the energy concentration factor as a combination of 
stress intensity factors.

\section{Energy release rate and energy concentration}
\indent

From relations (\ref{fvar}), (\ref{ji}), we deduce  that a good generalization 
of the $J$ integral of Rice (which is classically a  number) might a 
functional : 
$$\eta \ , \ supp \ \eta \subset \subset \Omega \ \mapsto \ - \ 
\frac{d}{dt}_{|t=0} \int_{\Omega} w(\nabla (\eu. \phi_{t}^{-1})) \mbox{ d}x 
\ \ , $$
where $\phi_{t}$ is the flow generated by $\eta$.

\begin{defi}
For any equilibrium state $\displaystyle (\eu,S) \in Eq(\eu_{0},K)$ and for any 
vector field $\eta$ which generates a one parameter flow $\displaystyle 
\phi_{t} \in \mathcal{D}(K)$, such that (there is a $T>0$ with) 
$\displaystyle S \subset \phi_{t}(S)$ for all $t \in [0,T]$, we define 
the energy release rate along the vector field $\eta$ by: 
\begin{equation}
ER(\eu,S)(\eta) \ = \ - \ \frac{d}{d t}_{|t=0} 
\ \int_{\Omega} w(\nabla (\eu \circ \phi_{t}^{-1}) ) \mbox{ d}x
\label{k2pre}
\end{equation}
\label{defer}
\end{defi}

Denote by $\mathcal{V}(K,S)$ the family of all vector fields $\eta$ 
generating a one parameter flow $\displaystyle 
\phi_{t} \in \mathcal{D}(K)$, such that there is a $T>0$ with  
$\displaystyle S \subset \phi_{t}(S)$ for all $t \in [0,T]$. Formally this set 
plays the role of the tangent space at the identity for the (infinite 
dimensional) semigroup of all $\phi \in \mathcal{D}(K)$ such that $S \subset 
\phi(S)$. 

Remark that $\displaystyle ER(\eu,S)(\eta)$ is  a linear expression in 
the variable $\eta$. Indeed, we have 
$$ER(\eu,S)(\eta)\ = \int_{\Omega} \left\{ \sigma(\nabla \eu)_{ij} \eu_{i,k} \eta_{k,j} \ - \ w(\nabla \eu) 
\ div \ \eta \right\} \mbox{ d}x \ \ . $$ 

 Nevertheless, the  set $\mathcal{V}(K,S)$ is not a 
vector space (mainly because the class of all $\phi \in \mathcal{D}(K)$ such 
that $S \subset \phi(S)$ is only a semigroup, and not a group). Therefore, the 
energy release rate is not a linear functional in a classical sense.

\begin{defi}
With the notations from definition \ref{defer}, the total variation of the 
energy release rate in  a open set $D \subset \Omega$  is defined by:
\begin{equation}
\mid ER \mid (\eu,S) (D) \ = \ \sup  \, ER(\eu,S)(\eta) \quad , 
\label{k2mes}
\end{equation}  
over all vector fields $\displaystyle \eta \in \mathcal{V}(K,S)$, with 
support in $D$, $supp \ \eta \subset D$, such that for all $x \in \Omega$ 
we have $\| \eta(x) \| \leq 1$. 
 
The function $\displaystyle \mid ER \mid (\eu,S)(\cdot)$ is positive and 
sub-additive, therefore induces  an outer measures over the $\sigma$-algebra of 
borelian sets in $\Omega$. 

We call this function the energy release rate 
associated to $\displaystyle (\eu,S) \in Eq(\eu_{0},K)$. 
\label{defeer}
\end{defi}

The number $\displaystyle \mid ER(\eu,S)\mid (D)$ measures the maximal 
elastic energy release rate that can be obtained by propagating the crack 
set $S$ inside the the set $D$, with sub-unitary speed, {\it by preserving it's
shape topologically}. 

In the case $n=2$, as explained in subsection \ref{neq2}, let $\displaystyle 
x_{0}$ be the crack tip of the crack set $S$, and $J$ the Rice integral. 
Then  for an open set 
$D \subset \Omega$ we have: 
\begin{enumerate}
\item[-] $\displaystyle \mid ER(\eu,S) \mid (D) = J$ if the crack tip belongs 
to $D$, that is $\displaystyle x_{0} \in D$, 
\item[-] $\displaystyle \mid ER(\eu,S) \mid (D) = 0$ if the crack tip does
not belong  to $D$. 
\end{enumerate}
For short, if we denote by $\displaystyle \delta x_{0}$ the Dirac measure 
centered at the crack tip $\displaystyle x_{0}$, we can write: 
$$\mid ER(\eu,S)\mid = J \ \delta x_{0} \quad .$$
It is therefore the
appropriate generalization of the Rice integral in three dimensions. 

Suppose that for any crack set $L$ and boundary displacement $\displaystyle 
\eu_{0}$ the surface energy has the expression: 
$$F(S;\eu_{0}) \ = \ G \hen(S) \quad .$$
Then $\displaystyle
CF(S,\eu_{0})(\Omega)$ is just $G$ times the perimeter (length if $n=3$) of 
the edge of the crack $S$ which is not contained in $K$ (technically, it is the Hausdorff measure
$\displaystyle \mathcal{H}^{n-2}$ of $dS \setminus K$). 

There is a mathematical formula which expresses the perimeter of the edge of 
an arbitrary  crack set $L$ as an "area release rate". Indeed, it is well known that 
the variation of the area of the crack set $\displaystyle \phi_{t}(L)$, along 
a one parameter flow generated by the vector field $\eta \in \mathcal{V}(K,L)$,   
has the expression: 
$$\frac{d}{d t}_{|t=0} \hen(\phi_{t}(S)) \ = \ \int_{S} div_{tan} \eta 
\mbox{ d}\hen(x) \quad , $$
where the operator $\displaystyle div_{tan}$ is the tangential divergence 
with respect to the surface $S$. If we denote by $\en$ the field of normals 
to the crack set $S$, then the expression of $\displaystyle div_{tan}$ 
operator is: 
$$div_{tan} \eta \ = \ \eta_{i,i} \ - \ \eta_{i,j} \en_{i} \en_{j} \quad . $$
Further, the perimeter of $dS \setminus K$, the edge of the crack set $S$
outside $K$, admits 
the following description, similar in principle to the expression of the 
elastic energy release rate given in definition \ref{defeer}: 
$$\mathcal{H}^{n-2}(dS\setminus K) \ = \ \sup \left\{ \int_{S} div_{tan} \eta 
\mbox{ d}\hen(x) \mbox{ : } \eta \in \mathcal{V}(K,S) , \, \forall x \in X \ \ 
\|\eta(x)\| \leq 1 \right\} \quad .$$
By putting together this expression of the perimeter, with relation (\ref{ji}), 
we obtain therefore the following proposition. 

\begin{prop}
If for any crack set $L$ we have $\displaystyle F(L;\eu_{0}) = G \hen(L)$ 
then for any absolute minimal state $\displaystyle (\eu,S) \in Absmin(\eu_{0}, 
K)$ such that $S \setminus K \not = \emptyset$ we have 
$$\mid ER(\eu,S) \mid(\Omega) \ = \ CF(\eu,S)(\Omega) \quad .$$
\label{pstrong}
\end{prop}

At this point let us remark that for a general equilibrium state in three 
dimensions 
$\displaystyle (\eu,S) \in Eq(\eu_{0},K)$  
 there is no obvious connection between the energy release rate 
 $\mid ER(\eu,S) \mid$, as in 
 definition \ref{defeer}, ant the elastic energy concentration $CE(\eu,S)$, 
 as in definition \ref{deficonc}.

The following theorem gives a relation between these two quantities. 

\begin{thm}
Let $\displaystyle (\eu,S) \in Eq(\eu_{0},K)$ be an equilibrium state 
of the brittle body with reference configuration $\Omega$, and $D \subset 
\Omega$ an arbitrary  open set. 
Then we have the following inequality: 
\begin{equation}
\mid ER(\eu,S) \mid (D) \ \leq \ CE(\eu,S)(D) \quad . 
\label{girv}
\end{equation}
\label{thmineq}
\end{thm}

\begin{rk}
For an arbitrary crack set $L$, we can't a priori deduce from the EH 
 the existence 
of a displacement $\eu'$ with  $\displaystyle (\eu', L) \in Adm(\eu_{0},K)$ and 
such that  for any other state  $\displaystyle (\ev, L) \in Adm(\eu_{0},K)$ we 
have 
$$\int_{\Omega} w(\nabla \eu') \mbox{ d}x  
\ \leq \int_{\Omega} w(\nabla \ev) \mbox{ d}x \quad .$$
From the mechanical point of view such an assumption is natural. There are 
mathematical results which supports this hypothesis, but as far as I know, 
not with the regularity needed in this paper. 
Fortunately, we shall not need to make such an assumption in order to prove 
theorem \ref{thmineq}. 
\end{rk}

\paragraph{Proof.}(First part) 
Let us consider an arbitrary  vector field $\eta \in \mathcal{V}(K,S)$, 
with compact support in $D$, such that for any $x \in \Omega$ we have 
$\|\eta(x)\| \leq 1$. 

In order to prove the theorem it is enough to show that  
\begin{equation}
ER(\eu,S)(\eta) \ \leq \ CE(\eu,S)(D) \quad .  
\label{girev}
\end{equation}
Indeed, suppose (\ref{girev}) is true for
any vector field $\eta \in \mathcal{V}(K,S)$, 
with compact support in $D$, such that for any $x \in \Omega$ we have 
$\|\eta(x)\| \leq 1$. Then, by  taking the supremum with respect to all 
such vector fields $\eta$ and using definition \ref{defeer},  we get the 
desired relation (\ref{girv}). 

The inequality (\ref{girev}) is a consequence of   
proposition \ref{pned3} and relation (\ref{tned5}),  which are of independent 
interest.  
We shall resume the proof of theorem \ref{thmineq}, by giving the proof of the 
inequality (\ref{girev}), after we  prove the before mentioned 
results. \quad $\square$

\vspace{.5cm}

Let $\displaystyle \phi_{t}$ be the one parameter flow generated 
by the vector field $\eta$. We can always find a curvilinear 
coordinate system $(\alpha_{1}, . . . , \alpha_{n-1}, \gamma)$ in the open 
set $D$  such that: 
\begin{enumerate}
\item[-] on the part of the edge $dS \cap \, supp \, \eta $ of the crack
set $S$ we have  $\gamma = 0$ ,  
\item[-] the surface $\gamma = t$ (constant) is the boundary of an open set 
$\displaystyle B_{t}$ such that 
$$\phi_{t}(S) \setminus S  \subset B_{t} \subset \, supp \, \eta \subset D \quad
,$$
\item[-] there exists $T>0$ such that for all $t \in [0,T]$ we have  
\begin{equation}
B_{t} \subset B(dS \cap D,t) \cap D \quad , 
\label{subuni}
\end{equation} 
\end{enumerate}
where $B(dS \cap D,t)$ is the tubular neighbourhood of $dS \cap D$, of radius 
$t$. 

Consider also the one parameter flow $\displaystyle \psi_{t}$, $t \geq 0$, 
which is equal to identity outside the open set $D$ and, in curvilinear
coordinates just introduced, it has the expression  
$$\psi_{t}(x(\alpha_{i},\gamma)) = x(\alpha_{i}, t+ \gamma) \quad .$$
Notice that $\displaystyle \psi_{t}(\Omega) = \Omega \setminus B_{t}$. 
We shall use these notations for proving that the elastic energy 
concentration is a kind of energy release rate, after the following result. 

\begin{prop}
With the notations made before, we have: 
\begin{equation}
\lim_{t \rightarrow 0} \frac{1}{t} \left( \int_{\Omega \setminus B_{t}} 
w(\nabla \eu) \mbox{
d}x - 
\int_{\Omega \setminus B_{t}} w(\nabla (\eu \circ \psi_{t}^{-1})) \mbox{ d}x 
\right) \ = 0 \quad .
\label{tned3}
\end{equation}
\label{pned3}
\end{prop}

\paragraph{Proof.} Recalling that $\displaystyle \psi_{t}(\Omega) = 
\Omega \setminus B_{t}$, we use the change of variables $\displaystyle 
x = \psi_{t}(y)$ to prove that  (\ref{tned3}) is equivalent with 
$$\lim_{t \rightarrow 0} \frac{1}{t} \left( \int_{\Omega} 
\left( w(\nabla\eu(y) (\nabla \psi_{t})^{-1}(y)) - w((\nabla
\eu)(\psi_{t}(y)) \right) \det \nabla \psi_{t}(y)  
\mbox{ d}y \right)  \ = 0 \quad .$$ 
The previous relation is just 
\begin{equation}
\frac{d}{d t}_{|t=0} \int_{\Omega} \left( w(\nabla\eu(y) (\nabla 
\psi_{t})^{-1}(y)) - w((\nabla\eu)(\psi_{t}(y)) \right) \det \nabla \psi_{t}(y)  
\mbox{ d}y \ = 0 \quad .
\label{tned4}
\end{equation}
We shall prove this from $\displaystyle (\eu,S) \in Eq(\eu_{0},K)$ and from 
an approximation argument. Notations from subsection \ref{configsect} will 
be in use. 

Denote by $\omega$ the vector field which generates the one parameter flow 
$\displaystyle \psi_{t}$. Let us
compute, using integration by parts: 
$$\frac{d}{d t}_{|t=0} \int_{\Omega} \left( w(\nabla\eu(y) (\nabla 
\psi_{t})^{-1}(y)) - w((\nabla\eu)(\psi_{t}(y))) \right) 
\det \nabla \psi_{t}(y)  
\mbox{ d}y \ = $$
\begin{equation}
= \  \int_{\Omega} \left( \sigma_{ij} \eu_{i,jk} \omega_{k} + 
\sigma_{ij} \eu_{i,k} \omega_{k,j} \right) \mbox{ d}y \quad .
\label{ned6}
\end{equation}

For any $\gamma > 0$, sufficiently small, choose 
a smooth scalar function $\displaystyle f^{\gamma}:\Omega \rightarrow
[0,1]$, such that:
\begin{enumerate}
\item[(a)] $\displaystyle f^{\gamma}(x) = 0$ for all $\displaystyle x \in
B_{\gamma}$, $\displaystyle f^{\gamma}(x) = 1$ for all $\displaystyle x \in
\Omega \setminus B_{2\gamma}$, 
\item[(b)] as $\gamma$ goes to $0$ we
have: 
$$\lim_{\gamma \rightarrow 0} \int_{\Omega} f^{\gamma} \left( \sigma_{ij} \eu_{i,jk} 
\omega_{k}+ \sigma_{ij} \eu_{i,k} \omega_{k,j} \right) \mbox{ d}y \ = \ 
\int_{\Omega} \left( \sigma_{ij} \eu_{i,jk} 
\omega_{k}+ \sigma_{ij} \eu_{i,k} \omega_{k,j} \right) \mbox{ d}y  \quad , $$
$$\lim_{\gamma \rightarrow 0} \int_{\Omega} f^{\gamma}_{,j} \sigma_{ij} \eu_{i,k} 
\omega_{k} \mbox{ d}y \ = 0 \quad . $$
\end{enumerate}
For all sufficiently small $\gamma >0$ it is true that: 
$$\int_{\Omega} \left( \sigma_{ij} \eu_{i,jk} \omega_{k}^{\gamma} + 
\sigma_{ij} \eu_{i,k} \omega^{\gamma}_{k,j} \right) \mbox{ d}y \ =$$
$$= \ \int_{\Omega} \left( f^{\gamma} \left( \sigma_{ij} \eu_{i,jk} 
\omega_{k}+ \sigma_{ij} \eu_{i,k} \omega_{k,j} \right)   + 
f^{\gamma}_{,j} \sigma_{ij} \eu_{i,k} \omega_{k}\right) \mbox{ d}y \quad .$$

Thus, from (a), (b) above  we get the equality: 
$$ \lim_{\gamma \rightarrow 0} \int_{\Omega} \left( \sigma_{ij} \eu_{i,jk} 
\omega_{k}^{\gamma}+ \sigma_{ij} \eu_{i,k} \omega_{k,j}^{\gamma} \right) \mbox{ d}y
 \ = \ \int_{\Omega} \left( \sigma_{ij} \eu_{i,jk} 
\omega_{k}+ \sigma_{ij} \eu_{i,k} \omega_{k,j} \right) \mbox{ d}y  \quad . $$
Recall that $(\eu,S)$ is an equilibrium state, therefore the stress field
$\sigma = \sigma(\nabla\eu)$ has divergence equal to $0$. 
Integration by parts shows that  for any sufficiently small  $\gamma > 0$ we 
have: 
$$\int_{\Omega} \left( \sigma_{ij} \eu_{i,jk} \omega_{k}^{\gamma} + 
\sigma_{ij} \eu_{i,k} \omega^{\gamma}_{k,j} \right) \mbox{ d}y \ = \ \int_{\Omega}- \sigma_{ij,j} \left(\eu_{i,k} \omega_{k}^{\gamma}\right) 
\mbox{ d} y \ = 0 \quad .$$ 
We obtained therefore the relation: 
$$\int_{\Omega} \left( \sigma_{ij} \eu_{i,jk} 
\omega_{k}+ \sigma_{ij} \eu_{i,k} \omega_{k,j} \right) \mbox{ d}y \ = 0 
\quad .$$
This is equivalent to relation (\ref{tned4}), by computation (\ref{ned6}).  \quad
$\square$

\vspace{.5cm}

A straightforward consequence of (\ref{tned3}) is that the elastic 
energy concentration 
is related to a kind of configurational energy release rate. Namely, we see 
that
$$\limsup_{t\rightarrow 0} \frac{1}{t}\int_{B_{t}} w(\nabla\eu) \mbox{ d}x \ =$$
\begin{equation}
= \ \limsup_{t \rightarrow 0} \frac{1}{t} \left( \int_{\Omega} 
w(\nabla \eu) \mbox{
d}x - 
\int_{\Omega \setminus B_{t}} w(\nabla (\eu \circ \psi_{t}^{-1})) \mbox{ d}x 
\right)  \quad .
\label{tned5}
\end{equation}

We turn back to the proof of theorem \ref{girv}. Recall that what it is left to 
prove is relation (\ref{girev}). 

\paragraph{Proof of (\ref{girev}).}
By construction, for all sufficiently small $t>0$  we have: 
$$\frac{1}{t} \int_{B(dS,t)\cap D} w(\nabla \eu) \mbox{ d}x \ \geq \ 
\frac{1}{t} \int_{B_{t}} w(\nabla \eu) \mbox{ d}x \quad .$$
because $\displaystyle B_{t} \subset B(dS, t) \cap D$. 
We write the right hand side member of this inequality as a sum of three  
terms: 
$$\frac{1}{t} \int_{B_{t}} w(\nabla \eu) \mbox{ d}x \ = $$
$$=\ \frac{1}{t} \left( \int_{\Omega} w(\nabla\eu) \mbox{ d}x -
\int_{\Omega} 
w(\nabla(\eu \circ \phi_{t}^{-1})) \mbox{ d}x \right) \ + $$
$$+\ \frac{1}{t} \left( \int_{\Omega} 
w(\nabla(\eu \circ \phi_{t}^{-1})) \mbox{ d}x - \int_{\Omega\setminus B_{t}} 
w(\nabla(\eu \circ \psi_{t}^{-1})) \mbox{ d}x \right) \ + $$
$$+\ \frac{1}{t} \left( \int_{\Omega\setminus B_{t}} 
w(\nabla(\eu \circ \psi_{t}^{-1})) \mbox{ d}x - \int_{\Omega\setminus B_{t}} 
w(\nabla(\eu)) \mbox{ d}x \right) \quad .$$
As $t$ goes to $0$, the first term converges to $EC(\eu,S)(\eta)$ and 
the third term converges to $0$ by proposition \ref{pned3}. We want to 
show that
\begin{equation}
\lim_{t \rightarrow 0} \frac{1}{t} \left( \int_{\Omega} 
w(\nabla(\eu \circ \phi_{t}^{-1})) \mbox{ d}x - \int_{\Omega\setminus B_{t}} 
w(\nabla(\eu \circ \psi_{t}^{-1})) \mbox{ d}x \right) \ = 0 \quad .
\label{tned7}
\end{equation}
The proof of this limit is  identical with the proof of proposition 
\ref{pned3}. Indeed, in that proof we worked with the one parameter flow 
$\displaystyle \psi_{t}$ generated by the vector field $\omega$. This one 
parameter flow is a semigroup (with respect to composition of functions), but 
after inspection of the proof it can be seen that we only used the following: 
 for any $x \in \Omega \setminus S$ 
$$\lim_{t \rightarrow 0} \psi_{t}(x) \ = \ x \quad \mbox{ and } \quad 
 \frac{d}{d t}_{|t=0} \psi_{t}(x) \ = \ \omega(x) \quad .$$
Therefore we can modify the proof of proposition \ref{pned3} by considering, 
instead of $\displaystyle \psi_{t}$, the diffeomorphisms $\displaystyle
\lambda_{t}$ defined by:
$$\lambda_{t} \ = \ \psi_{t}\circ \phi_{t}^{-1} \quad .$$
The rest of the proof goes exactly as before, thus leading us to relation 
(\ref{tned7}). 

Eventually, we have: 
$$CE(\eu,S)(D) \ = \ \limsup_{t \rightarrow 0} \frac{1}{t} \int_{B(dS,t)\cap D}
w(\nabla \eu) \mbox{ d}x \ \geq$$
$$\geq \ \limsup_{t \rightarrow 0} \frac{1}{t} \int_{B(dS,t)\cap D}
w(\nabla \eu) \mbox{ d}x \ = \ ES(\eu,S)(\eta) \  + $$
$$+ \ \lim_{t \rightarrow 0} \frac{1}{t} \left( \int_{\Omega} 
w(\nabla(\eu \circ \phi_{t}^{-1})) \mbox{ d}x - \int_{\Omega\setminus B_{t}} 
w(\nabla(\eu \circ \psi_{t}^{-1})) \mbox{ d}x \right) \  + \ $$
$$+ \ \lim_{t \rightarrow 0} \frac{1}{t} \left( \int_{\Omega\setminus B_{t}} 
w(\nabla(\eu \circ \psi_{t}^{-1})) \mbox{ d}x - \int_{\Omega\setminus B_{t}} 
w(\nabla(\eu)) \mbox{ d}x \right) \  = \ ES(\eu,S)(\eta)$$
and (\ref{girev}) is therefore proven. \quad $\square$

\section{A constraint on some minimal solutions}

Let us consider now a solution of the model of brittle crack propagation
described in section \ref{themodels}. More precisely, for given boundary 
conditions $ \displaystyle \eu_{0}(t)$ and initial crack set $K$, 
we shall call a solution 
$\displaystyle (\eu(t), S_{t}) \in Eq(\eu_{0}(t),S_{t})$ of the model described 
by axioms (A1),..., (A5), by the name "equilibrium solution". Likewise,  
a solution 
$\displaystyle (\eu(t), S_{t}) \in Absmin(\eu_{0}(t),S_{t})$ of the model 
described by axioms (A1),(A2),(A3'),(A4), will be called a "minimal solution".

We shall  deal with a minimal solution $\displaystyle (\eu(t), S_{t}) \in 
Absmin(\eu_{0}(t),S_{t})$ for which the crack set $\displaystyle S_{t}$
propagates smoothly, {\it without topological changes}. 
Namely we shall suppose that there exists a vector field 
$\eta$ with compact support in $\Omega$, such that for all $t \in [0,T]$ we 
have $\displaystyle S_{t} = \phi_{t}(K)$, where $\displaystyle \phi_{t}$ is 
the one parameter flow generated by $\eta$. 

Because the problem is quasistatic, time enters only as a parameter, therefore 
we may suppose moreover that for all $x \in \Omega$ we have $\eta(x) \leq 1$. 

At each moment $t \in [0,T]$ we shall have $\displaystyle \eta \circ \phi_{t} \in \mathcal{V}(K, S_{t})$.

\begin{thm}
Suppose that for any crack set $L$ and boundary displacement $\displaystyle 
\eu_{0}$ the surface energy has the expression: 
$$F(S;\eu_{0}) \ = \ G \hen(S) \quad .$$
Let $\displaystyle (\eu(t), S_{t}) \in 
Absmin(\eu_{0}(t),S_{t})$ be a minimal solution, with $\displaystyle S_{0}=K$, 
 such that exists a vector field $\eta$ with $\|\eta(x)\| \leq 1$ for all 
 $x \in \Omega$ and for all $t \in [0,T]$ we have $\displaystyle S_{t} = 
 \phi_{t}(K)$, where $\displaystyle \phi_{t}$ is the one parameter flow 
 generated by $\eta$. 
 
 Then for any $t \in [0,T]$  and any open set $D \subset \Omega$ 
 we have the equalities: 
 $$\mid ER(\eu(t),\phi_{t}(S)) \mid (D) \ = \  EC(\eu(t),\phi_{t}(S))(D) \ =$$
 \begin{equation}
 = \ CF(\phi_{t}(S);\eu_{0}(t))(D) \ = \ G \mathcal{H}^{n-2}(dS \setminus K) \quad
.
\label{maineq}
\end{equation}
\label{mainthm}
\end{thm}

\paragraph{Proof.}
Theorems \ref{thmineq} and \ref{firsthm} tell us that for any open set 
$D \subset \Omega$, and for any $t \in [0,T]$ we have 
$$\mid ER(\eu(t),\phi_{t}(S)) \mid (D) \leq EC(\eu(t),\phi_{t}(S))(D) \leq 
CF(\phi_{t}(S);\eu_{0}(t))(D) \quad .$$
Proposition \ref{pstrong} tells that 
$$CF(\phi_{t}(S); \eu_{0}(t))(\Omega) \ = \  \mid ER(\eu(t), \phi_{t}(S)) \mid 
(\Omega) \quad .$$
We deduce that for any open set 
$D \subset \Omega$, and for any $t \in [0,T]$ the string of equalities 
(\ref{maineq}) is true. \quad $\square$

This result is natural in two dimensional linear elasticity. Nevertheless, in 
the case of three dimensional elasticity, the constraints on the elastic energy
concentration provided by theorem \ref{mainthm} might be too hard to satisfy. 

Indeed, from (\ref{maineq}) we deduce that in  particular the elastic energy 
concentration has to be absolutely continuous with respect to the perimeter 
measure of the edge of the crack.

\section{Conclusions}

We have proposed a general model of brittle crack propagation based on 
Mumford-Shah functionals. We have defined equilibrium and absolute minimal
solutions of the model.

By a combination of analytical and configurational analysis, we defined 
measures of energy release rate and energy concentrations for equilibrium and 
absolute minimal solutions and we have shown that there is a difference between 
such solutions, as shown mainly by theorems \ref{firsthm}, \ref{thmineq} and
\ref{mainthm}.


\begin{thebibliography}{10}
\bibitem{A1}  L. Ambrosio, Variational problems in SBV and image segmentation, Acta Appl.  
Mathematic\ae \quad  17, 1989,1-40

\bibitem{A2}  L. Ambrosio, Existence Theory for a New Class of Variational Problems, 
Arch. Rational Mech. Anal., vol. 111, 1990, 291-322



\bibitem{BCDM}  G. Bellettini, A. Coscia, G. Dal Maso, Compactness and lower 
semicontinuity properties in $\eSBD(\Omega)$, {\it Mathematische Zeitschrift}, 
vol. 228, 2, 337-351, 1998



\bibitem{bu} M. Buliga, Variational formulations in brittle fracture mechanics, 
Ph.D. Thesis, Institute of Mathematics of the Romanian Academy, 1997

\bibitem{bu3} M. Buliga, Energy concentration and brittle crack propagation, 
J. of Elasticity, {\bf 52}, 3, 201-238, 1999

\bibitem{bu2} M. Buliga, Brittle crack propagation based on an optimal energy 
balance, {\it Rev. Roum. des Math. Pures et Appl.}, 45, 2, 201-209, 2001





\bibitem{dft05} G. Dal Maso, G. Francfort, R. Toader, Quasistatic crack growth 
in nonlinear elasticity, {\it Arch. Rat. Mech. Anal.}, 176, 165-225, 2005

\bibitem{DGA}  E. De Giorgi, L. Ambrosio, Un nuovo funzionale del calcolo delle variazioni, Atti Accad. Naz. 
Lincei Rend. Cl. Sci. Fis. Mat. Natur., 82, 1988, 199-210





\bibitem{FMa} G. Francfort, J.-J. Marigo, Stable damage evolution in a brittle continuous medium, 
Eur. J. Mech., A/Solids, 12, 149-189, 1993

\bibitem{FMa2} G. Francfort, J.-J. Marigo, Revisiting brittle fracture as an 
energy minimization problem, J. Mech. Phys. Solids, 46, 1319-1342, 1998



\bibitem{mielke} A. Mielke, Evolution in rate-independent systems (Ch. 6). In 
C. Dafermos, E. Feireisl, eds., {\it Handbook of Differential Equations, 
Evolutionary Equations, vol. 2}, 461-559, Elsevier B.V., Amsterdam, 2005



\bibitem{MS} D. Mumford, J. Shah, Optimal approximation by piecewise smooth functions and associated 
variational problems, Comm. on Pure and Appl. Math., vol. XLII, 
no. 5, 1989


\end{thebibliography}
\end{document}